\documentclass[11pt]{article}
\usepackage{latexsym,amssymb}
\def\demo{\noindent{\bf Proof. }}
\def\QED{\hfill$\Box$}

\newtheorem{Theorem}{Theorem}[section]
\newtheorem{Lemma}[Theorem]{Lemma}
\newtheorem{Corollary}[Theorem]{Corollary}
\newtheorem{Proposition}[Theorem]{Proposition}
\newtheorem{Remark}[Theorem]{Remark}
\newtheorem{Example}[Theorem]{Example}

\newtheorem{Definition}[Theorem]{Definition}

\begin{document}
\topmargin3mm
\hoffset=-1cm
\voffset=-1.5cm
\

\medskip

\begin{center}
{\large\bf Ring graphs and complete intersection toric ideals}
\vspace{6mm}\\

\medskip

Isidoro Gitler, Enrique Reyes
\footnote{Partially supported by COFAA-IPN.}, and Rafael H. Villarreal 
\footnote{This work was partially supported
by CONACyT grants 49251-F, 49835-F, and SNI.}
\\ 
{\small Departamento de Matem\'aticas}\vspace{-1mm}\\ 
{\small Centro de Investigaci\'on y de Estudios Avanzados del
IPN}\vspace{-1mm}\\   
{\small Apartado Postal 14--740}\vspace{-1mm}\\ 
{\small 07000 M\'exico City, D.F.}\vspace{-1mm}\\ 
{\small e-mail: {\{{\tt igitler,ereyes,vila}\}{\tt @math.cinvestav.mx}}}\vspace{4mm}
\end{center}
\date{}

\begin{abstract} 
\noindent We study the family of graphs whose number of
primitive cycles equals its cycle rank. It is shown that this family
is precisely the family of ring graphs. Then we study
the complete intersection property of toric ideals of 
bipartite graphs and oriented graphs. An interesting application is
that complete intersection toric 
ideals of bipartite graphs correspond to ring graphs and that these
ideals are minimally generated by Gr\"obner bases. We prove that any
graph can be oriented such that its toric ideal is a complete
intersection with a universal Gr\"obner basis 
determined by the cycles. It turns out that bipartite ring graphs are exactly the
bipartite graphs that have 
complete intersection toric ideals for any orientation.  
\end{abstract}

\section{Introduction} 
Let $G$ be a graph (no loops or multiple edges) with $n$ vertices and $q$ edges, 
and let ${\rm frank}(G)$ be the number of primitive cycles of $G$,
i.e., cycles without chords. The number ${\rm frank}(G)$ is 
called the {\it free rank\/} of $G$ and the number 
${\rm rank}(G)=q-n+r$ is called the {\it cycle rank\/} of $G$, where
$r$ is the number of connected components of $G$. The cycle rank of
$G$ can be expressed as the dimension of the cycle space
of $G$. These two numbers 
satisfy ${\rm rank}(G)\leq {\rm frank}(G)$, as is seen in 
Proposition~\ref{rank-leq-frank}. The aim of this paper is to study and
classify the family of graphs where the equality occurs. It will
turn out that this family is precisely the family of ring graphs. The
precise definition of a ring graph can be found in
Section~\ref{ring-graphs}. Roughly speaking {\it ring graphs\/} 
can be obtained starting with a cycle and subsequently attaching
paths of  length at least two that meet graphs already constructed in
two adjacent vertices.   

The contents of this paper are as follows. Before stating our main
results, recall that a graph $G$ has the {\it primitive cycle
property\/} (PCP) if any two primitive cycles intersect in at most one
edge. A {\it subdivision} of a graph is any graph that can be
obtained from the original graph by replacing edges
by paths.  As usual we denote the complete graph on $n$ vertices 
by ${\cal K}_n$. In Section~\ref{ring-graphs}, which is the core of the paper,
we prove the following implications for any graph $G$: 
$$
\begin{array}{ccccccc}
\mbox{outerplanar}&\Rightarrow &\mbox{ring graph}&
\Leftrightarrow & \mbox{PCP}+\mbox{contains no }& &\\ 
&&\Updownarrow& &\mbox{subdivision of }\mathcal{K}_{4} &\Rightarrow& \mbox{planar}\\ 
&&{\rm rank}={\rm frank}& &\mbox{as a subgraph} & 
\end{array}
$$ 

These purely graph theoretical results are applied in 
Sections~\ref{simplegraphs} and \ref{orientedgraphs}, 
where graphs with complete intersection toric
ideals are studied, both in the oriented and unoriented case.  
For bipartite graphs the equality ${\rm rank}(G)={\rm frank}(G)$ is
related to these special types of toric ideals as we explain below. 

Let $R=k[x_1,\ldots,x_n]$ be a polynomial ring 
over a field $k$ and let $G$ be a graph with vertex set 
$V(G)=\{x_1,\ldots,x_n\}$ and edge set $E(G)=\{t_1,\ldots,t_q\}$.  
The {\it edge subring\/} of $G$ is the $k$-subalgebra of $R$:
\[
k[G]=k[\{x_ix_j\vert\, x_i\, \mbox{ is adjacent to }x_j\}]
\subset R.
\]
There is an epimorphism of $k$-algebras
\[ 
\varphi\colon k[t_1,\ldots,t_q] \longrightarrow k[G], \ \ \ 
\{x,y\}\longmapsto xy,
 \]
where $k[t_1,\ldots,t_q]$ is a polynomial ring. The kernel of $\varphi$, denoted by
$P(G)$, is called the {\it toric ideal\/} 
of $G$. Toric ideals of graphs are studied in Section~\ref{simplegraphs}. 
The height of $P(G)$ is equal to $g=q-{\rm rank}(A_G)$, 
where $A_G$ is the incidence matrix of $G$. By a 
result of Krull \cite{AM} the ideal $P(G)$ cannot be generated by
less than $g$
polynomials. The toric ideal of $G$ is
called a {\it complete intersection\/} if it can be generated 
by $g$ polynomials. The complete
intersection property of $P(G)$ was first studied in 
\cite{luisa-tor,aron-jac}, and later in \cite{accota-gv,katzman}. 

An interesting result of Simis \cite{aron-jac}
shows that if $G$ is a bipartite graph, then 
${\rm rank}(G)={\rm frank}(G)$ if and only if $P(G)$ is 
a complete intersection. Thus by describing the graphs where equality
occurs, we are in particular describing the toric ideals of bipartite
graphs that are complete intersections 
(see Corollary~\ref{aug20-05}). We prove that 
complete intersection toric ideals of $2$-connected bipartite graphs 
are minimally generated by Gr\"obner bases (see
Corollary~\ref{jul16-05}).

In Section~\ref{orientedgraphs} we
introduce and study toric ideals of oriented graphs and their
Gr\"obner bases. To the best
of our knowledge these toric ideals have not been studied much except
for the case of acyclic tournaments \cite{ishizeki}. 
Oriented graphs share some properties with bipartite
graphs. For instance in both cases their incidence matrices are 
totally unimodular. This is a key fact to understand 
the Gr\"obner bases of toric ideals of oriented graphs 
(see Lemma~\ref{jul17-05}). We prove that 
the toric ideal of any oriented graph is completely determined 
by  its primitive cycles and has a 
universal Gr\"obner basis determined by the cycles 
(see Proposition~\ref{grobnbas} and Corollary~\ref{jul22-1-05}). It 
is shown that toric ideals of oriented ring graphs are complete 
intersections for any
orientation. As an interesting consequence 
of the results of Section~\ref{ring-graphs} we obtain that for 
bipartite graphs this property characterizes ring graphs
(see Corollary~\ref{aug10-05}). One of our main results shows that
any graph has an acyclic orientation such that the corresponding 
toric ideal is a complete intersection (see Theorem~\ref{icd}). 

The paper is essentially self contained. For unexplained terminology
and notation on graph theory we refer to \cite{diestel,Har}. Our main
references for edge subrings are \cite{Vi3,monalg}.

\section{Ring graphs}\label{ring-graphs}
Let $G$ be a graph with $n$ vertices and $q$ edges. We denote the 
vertex set and edge set of $G$ by $V(G)=\{x_1,\ldots,x_n\}$ and 
$E(G)=\{t_1,\ldots,t_q\}$ respectively. Recall that a $0$-chain 
(resp. $1$-chain) of $G$ is a formal linear combination 
$\sum a_ix_i\ \  (\mbox{resp.}\ \sum b_it_i$) of vertices (resp.
edges), where $a_i\in \mathbb{Z}_2$ (resp. $b_i\in \mathbb{Z}_2$). 
The {\it boundary operator\/} is the 
linear map $\partial\colon C_1\rightarrow C_0$ defined by 
$$\partial
(\{x,y\})=x+y,
$$
where $C_i$ is the $\mathbb{Z}_2$-vector space of $i$-chains. A {\it cycle
vector\/} is a $1$-chain of the form  
$t_1+\cdots+t_r$ where $t_1,\ldots,t_r$ are the edges of a cycle 
of $G$. 
The 
{\it cycle space\/} ${\cal Z}({G})$
of $G$ over $\mathbb{Z}_2$ is equal to ${\rm ker}(\partial)$. The 
vectors in ${\cal Z}({G})$ can be regarded as a set of
edge-disjoint cycles. A {\it cycle basis\/}
for $G$ is a basis
for ${\cal Z}({G})$ which consists entirely of cycle vectors, such a
basis can be constructed as follows:  

\begin{Remark}{\rm \cite[pp. 38-39]{Har}}\label{spanntreee}\rm\ If
$G$ is connected, then $G$ has a spanning tree $T$. 
The subgraph of $G$ consisting of $T$ and any edge 
in $G$ not in $T$ has exactly one cycle, the collection of all cycle 
vectors of cycles
obtained in this way form a cycle basis for $G$. Hence 
${\rm dim}_{\mathbb{Z}_2}\, {\cal Z}({G}) =q-n+r$ if $G$ is a graph
with $r$ connected components. 
\end{Remark}

Let $c$ be a cycle of $G$. A {\it chord\/} of $c$ is any edge of $G$ joining two 
non adjacent vertices of $c$. 
A cycle without chords is called 
{\it primitive\/}. The number 
${\rm dim}_{\mathbb{Z}_2}\, {\cal Z}({G})$ is called the {\it cycle
rank} of $G$ and is denoted by ${\rm
rank}(G)$. The number of primitive cycles of a graph $G$, denoted by 
${\rm frank}(G)$, is called the {\it free rank} of $G$.  

\begin{Proposition}\label{rank-leq-frank} If $G$ is a graph, 
then ${\cal Z}(G)$ is generated by cycle vectors of 
primitive cycles. In 
particular ${\rm rank}(G)\leq {\rm frank}(G)$.
\end{Proposition}

\demo Let $\mathbf{c}_1,\ldots,\mathbf{c}_r$ be a cycle basis for the
cycle space of 
$G$ and let $c_1,\ldots,c_r$ be the corresponding cycles of $G$. 
It suffices to notice that if 
some $c_j$ has a chord, we can write
$\mathbf{c}_j=\mathbf{c}_j'+\mathbf{c}_j''$, where $\mathbf{c}_j'$
and $\mathbf{c}_j''$ are cycle vectors of cycles 
of length smaller than that of $c_j$. \QED

\begin{Corollary}\label{vila-isi} Let $G$ be a graph. Then the following are
equivalent\/{\rm :}\vspace{-2mm}
\begin{description}
\item{\rm (a)} ${\rm rank}(G)={\rm frank}(G)$. \vspace{-2mm}
\item{\rm (b)} The set of cycle vectors of primitive cycles is a
basis for 
${\cal Z}(G)$. \vspace{-2mm}
\item{\rm (c)} The set of cycle vectors of primitive cycles is linearly independent.
\end{description}
\end{Corollary}

\demo (a) $\Rightarrow$ (b): By Proposition~\ref{rank-leq-frank} there
is a basis $\mathcal{B}$ of $\mathcal{Z}(G)$
consisting of cycle vectors of primitive cycles. By hypothesis 
${\rm rank}(G)={\rm frank}(G)$. Thus $\mathcal{B}$ is the set of all 
cycle vectors of primitive cycles and $\mathcal{B}$ is a basis. That
(b) implies (c) and (c) implies (a) are also very easy to prove. \QED

\medskip

Let $G$ be a graph. A vertex $v$ (resp. an edge $e$) of $G$  is called a 
{\it cutvertex\/} (resp. {\it
bridge\/}) if the number of connected components
of $G\setminus\{v\}$ (resp. $G\setminus\{e\}$) is larger than that of
$G$. A 
maximal connected subgraph of $G$ without cutvertices 
is called a {\it block\/}. A graph $G$ 
is $2$-{\it connected\/} if $|V(G)|>2$ and $G$ has no 
cutvertices. Thus a block of $G$ is either a maximal 
$2$-connected subgraph, a bridge or an isolated vertex. By their maximality, 
different blocks of $G$ intersect in at most one vertex, 
which is then a cutvertex of $G$. Therefore every 
edge of $G$ lies in a unique block, and $G$ is the 
union of its blocks. 

\begin{Lemma}\label{mar22-05} Let $G$ be a graph and let
$G_1,\ldots,G_r$ be its blocks. Then ${\rm rank}(G)={\rm frank}(G)$ if and
only if  ${\rm rank}(G_i)={\rm frank}(G_i)$ for all $i$. 
\end{Lemma}

\demo $\Rightarrow$) Let $G_i$ be any block of $G$. We may assume 
$|V(G_i)|>2$, otherwise ${\rm rank}(G_i)={\rm frank}(G_i)=0$. If $c$
is a primitive cycle of $G_i$, then by the maximality condition of 
a block one has that $c$ is also a primitive cycle of $G$. Thus by 
Corollary~\ref{vila-isi} the set of cycle vectors of primitive cycles
of $G_i$ is linearly independent and ${\rm rank}(G_i)={\rm frank}(G_i)$. 
 
$\Leftarrow$) Let $\mathcal{B}_i$ and $\mathcal{B}$ be the set of
cycle vector of primitive cycles of $G_i$ and $G$ respectively. 
As $\cup_{i=1}^r\mathcal{B}_i$ is linearly independent, by
Corollary~\ref{vila-isi} it suffices to
prove that $\cup_{i=1}^r\mathcal{B}_i=\mathcal{B}$. In the first part
of the proof we have already observed that 
$\cup_{i=1}^r\mathcal{B}_i\subset \mathcal{B}$. To prove the equality take
any cycle vector $\mathbf{c}$ of a primitive cycle $c$ of $G$. Since
$c$ is a $2$-connected subgraph, it must be contained in some block
of $G$, i.e., in some $G_i$. Thus $c$ is a primitive cycle 
of $G_i$, so $\mathbf{c}$ is in $\mathcal{B}_i$. \QED

\begin{Definition}\rm Given a graph $H$, we call a path $\cal P$ an
$H$-{\it path\/} if $\cal P$ is non-trivial and meets
$H$ exactly in its ends.
\end{Definition}

In order to describe, in graph theoretical terms, 
the family of graphs satisfying the equality ${\rm rank}(G)={\rm
frank}(G)$ we need to introduce another notion.


\begin{Definition}\rm A graph $G$ is a {\it ring graph\/}
 if each block of $G$ which is not a bridge or a vertex can be constructed from a 
cycle by successively adding $H$-paths of length at least $2$ that
meet graphs $H$ already 
constructed in two adjacent vertices.
\end{Definition}

Families of ring graphs include forests and cycles. These graphs are
planar by construction. 


\begin{Remark}\label{march23-love}\rm Let $G$ be a $2$-connected ring
graph and let $c$  
be a fixed primitive cycle of $G$, then $G$ can be constructed from $c$ by 
successively adding $H$-paths of length at least $2$ that
meet graphs $H$ already constructed in two adjacent vertices.
\end{Remark}  

A graph $H$ is called a {\it subdivision}
 of a graph $G$ if $H=G$ or $H$ arises from $G$ by replacing
edges by paths. 

\begin{Lemma}{\rm\cite[Lemma~7.78, p.~387]{aigner}}\label{dirac05}
Let $G$ be 
a graph with vertex set $V\!$.
If $G$ is $2$-connected and $\deg(v)\geq 3$ for all $v\in V\!$, then 
$G$ contains a subdivision of ${\cal K}_4$ as a subgraph.
\end{Lemma}

\begin{Lemma}\label{lina} Let $G$ be a graph. If ${\rm rank}(G)={\rm
frank}(G)$ and $x,y$ are two non adjacent vertices of $G$, then there
are at most two vertex disjoint paths joining $x$ and $y$.
\end{Lemma}

\demo Assume that there are three vertex disjoint paths joining $x$ and $y$:
\begin{eqnarray*}
{\cal P}_1=\{x,x_1,\ldots,x_r,y\},&
&{\cal P}_2=\{x,z_1,\ldots,z_t,y\},\\
{\cal P}_3=\{x,y_1,\ldots,y_s,y\},& &
\end{eqnarray*}
where $r,s,t$ are greater or equal than $1$. We may assume 
that the sum of the lengths of the ${\cal P}_i$'s is minimal. Consider the cycles
\begin{eqnarray*}
c_1&=&\{x,x_1,\ldots,x_r,y,z_t,\ldots,z_1,x\},\\
c_2&=&\{x,z_1,\ldots,z_t,y,y_s,\ldots,y_1,x\},\\
c_3&=&\{x,x_1,\ldots,x_r,y,y_s,\ldots,y_1,x\}.
\end{eqnarray*}
Thus we are in the following situation:

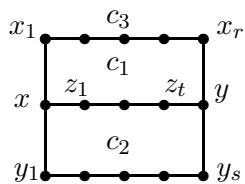
\setlength{\unitlength}{.035cm}
\thicklines
\begin{picture}(80,100)(-150,-25)
\put(-30,22){\circle*{4.1}}
\put(-30,22){\line(1,0){59}}
\put(-42,21){$x$}
\put(-30,47){\circle*{4.1}}
\put(-23,27){$z_1$}
\put(-30,21){\line(0,-1){26}}
\put(-30,21){\line(0,1){26}}
\put(15,27){$z_t$}
\put(0,22){\circle*{4.1}}
\put(-15,22){\circle*{4.1}}
\put(15,22){\circle*{4.1}}
\put(30,22){\circle*{4.1}}
\put(34,25){$y$}
\put(0,47){\circle*{4.1}}
\put(-30,47){\circle*{4.1}}
\put(-15,47){\circle*{4.1}}
\put(15,47){\circle*{4.1}}
\put(30,47){\circle*{4.1}}
\put(30,47){\circle*{4.1}}
\put(-30,47){\line(1,0){59}}

\put(-30,-5){\circle*{4.1}}
\put(-42,-5){$y_1$}
\put(-30,-5){\line(1,0){59}}

\put(-44,49){$x_1$}
\put(35,49){$x_r$}

\put(0,-5){\circle*{4.1}}
\put(-15,-5){\circle*{4.1}}
\put(15,-5){\circle*{4.1}}
\put(30,-5){\circle*{4.1}}
\put(30,-5){\line(0,1){52}}
\put(35,-5){$y_s$}

\put(-7,35){$c_1$}

\put(-7,5){$c_2$}

\put(-7,53){$c_3$}

\end{picture}

\noindent Observe that, by the choice of the ${\cal P}_i$'s, a chord
of the cycle $c_1$ (resp. $c_2$, $c_3$) must join $x_i$ and $z_j$
(resp. $z_i$ and $y_j$, $x_i$ and $y_j$) for some $i,j$. If $c_1$ is
not primitive, we can write 
$$
\mathbf{c}_1=\mathbf{a}_1+\cdots+\mathbf{a}_{n_1}
$$
for some distinct cycle vectors 
$\mathbf{a}_1,\ldots,\mathbf{a}_{n_1}$ of primitive cycles 
$\mathrm{a}_1,\ldots,\mathrm{a}_{n_1}$
such that each cycle $\mathrm{a}_i$ contains at least one edge of the 
form $\{x_j,z_k\}$. Similarly if $c_2$ (resp. $c_3$) is not primitive
we can write:
$$
\mathbf{c}_2=\mathbf{b}_1+\cdots+\mathbf{b}_{n_2}\ \mbox{ (resp. }\ 
\mathbf{c}_3=\mathbf{d}_1+\cdots+\mathbf{d}_{n_3})
$$
for some distinct cycle vectors
$\mathbf{b}_1,\ldots,\mathbf{b}_{n_2}$ (resp.
$\mathbf{d}_1,\ldots,\mathbf{d}_{n_3}$) of primitive cycles 
such that each cycle $\mathrm{b}_i$ (resp. $\mathrm{d}_i$) contains at
least one edge of the  
form $\{z_j,y_k\}$ (resp. $\{x_j,y_k\}$). Therefore we
can write 
$$
\mathbf{c}_1=\sum_{i=1}^{n_1}\mathbf{a}_i,\ \ 
\mathbf{c}_2=\sum_{i=1}^{n_2}\mathbf{b}_i,\ \ 
\mathbf{c}_3=\sum_{i=1}^{n_3}\mathbf{d}_i
$$
where $\mathbf{a}_1,\ldots,\mathbf{a}_{n_1},
\mathbf{b}_1,\ldots,\mathbf{b}_{n_2},\mathbf{d}_1,\ldots,\mathbf{d}_{n_3}$
are distinct cycle vectors of primitive cycles of $G$. Thus 
from the equality $\mathbf{c}_3=\mathbf{c}_1+\mathbf{c}_2$ we get a
non trivial linear relation of the set $\mathcal{B}$ of 
cycle vectors of primitive cycles, i.e., $\mathcal{B}$ is linearly
dependent, a contradiction to Corollary~\ref{vila-isi}. \QED

\begin{Lemma}\label{march11-03} Let $G$ be a graph. If ${\rm rank}(G)={\rm
frank}(G)$, then $G$ has the primitive cycle property. 
\end{Lemma}

\demo Let $c_1,c_2$ be two distinct primitive cycles. Assume that
$c_1$ and $c_2$ intersect in at least two edges. Thus $c_1$ and $c_2$ 
must intersect in at least two non adjacent vertices $u,v$. The cycle
$c_2$ can be written as:
$$ 
c_2=\{u=u_0,u_1,\ldots,u_s,v=u_{s+1},v_1,\ldots,v_m,u\}.
$$
At least one of the paths $\mathcal{P}_1=\{u,u_1,\ldots,u_s,v\}$, 
$\mathcal{P}_2=\{v,v_1,\ldots,v_m,u\}$ that form the cycle $c_2$ 
must contain a vertex not in
$c_1$, otherwise $c_1=c_2$. Assume that the path $\mathcal{P}_1$ has
this property. Hence there is $u_k\notin c_1$ such that $u_i\in c_1$ for 
$i<k$, and there is $u_\ell\in c_1$, with $k<\ell$, 
such that $u_i\notin c_1$ 
for $k\leq i<\ell$. Hence there are two non
adjacent vertices $x=u_{k-1},y=u_\ell$ in $c_1$ and a 
path ${\cal P}=\{x,u_k,\ldots,u_{\ell-1},y\}$ 
of length at least two that intersect $c_1$ in exactly 
the vertices $x,y$:

\setlength{\unitlength}{.035cm}
\thicklines
\begin{picture}(80,100)(-150,-25)
\put(-30,22){\circle*{4.1}}
\put(-30,22){\line(1,0){59}}
\put(-42,21){$x$}
\put(-30,47){\circle*{4.1}}
\put(-30,21){\line(0,-1){26}}
\put(-30,21){\line(0,1){26}}
\put(0,22){\circle*{4.1}}
\put(-15,22){\circle*{4.1}}
\put(15,22){\circle*{4.1}}
\put(30,22){\circle*{4.1}}
\put(34,25){$y$}
\put(0,47){\circle*{4.1}}
\put(-30,47){\circle*{4.1}}
\put(-15,47){\circle*{4.1}}
\put(15,47){\circle*{4.1}}
\put(30,47){\circle*{4.1}}
\put(30,47){\circle*{4.1}}
\put(-30,47){\line(1,0){59}}

\put(-30,-5){\circle*{4.1}}
\put(-30,-5){\line(1,0){59}}

\put(-44,49){$u_k$}
\put(35,49){$u_{\ell-1}$}

\put(0,-5){\circle*{4.1}}
\put(-15,-5){\circle*{4.1}}
\put(15,-5){\circle*{4.1}}
\put(30,-5){\circle*{4.1}}
\put(30,-5){\line(0,1){52}}
\put(-7,5){$c_1$}
\put(-7,53){$\cal P$}
\end{picture}

\noindent This contradicts Lemma~\ref{lina}. 
\QED

\medskip

\begin{Lemma}\label{addlemma2} Let $G$ be a graph. If $G$ satisfies
{\rm PCP} and $G$ 
does not contain a subdivision of ${\cal K}_4$ as a subgraph, then for
any two non adjacent vertices $x,y$ of $G$ there 
are at most two vertex disjoint paths joining $x$ and $y$.
\end{Lemma}

\demo Assume that there are three vertex disjoint paths joining $x$ and $y$:
$$
{\cal P}_1=\{x,x_1,\ldots,x_r,y\},\  {\cal
P}_2=\{x,z_1,\ldots,z_t,y\},\ {\cal P}_3=\{x,y_1,\ldots,y_s,y\},
$$
where $r,s,t$ are greater or equal than $1$. We may assume 
that the sum of the lengths of the ${\cal P}_i$'s is minimal. Consider the cycles
\begin{eqnarray*}
c_1=\{x,x_1,\ldots,x_r,y,z_t,\ldots,z_1,x\},& c_2=
\{x,z_1,\ldots,z_t,y,y_s,\ldots,y_1,x\},& \\
c_3=\{x,x_1,\ldots,x_r,y,y_s,\ldots,y_1,x\}.& &
\end{eqnarray*}
Thus we are in the following situation:

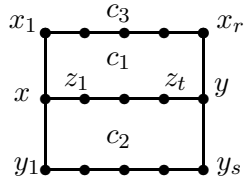
\begin{picture}(80,100)(-150,-25)
\setlength{\unitlength}{.035cm}
\thicklines
\put(-30,22){\circle*{4.1}}
\put(-30,22){\line(1,0){59}}
\put(-42,21){$x$}
\put(-30,47){\circle*{4.1}}
\put(-23,27){$z_1$}
\put(-30,21){\line(0,-1){26}}
\put(-30,21){\line(0,1){26}}
\put(15,27){$z_t$}
\put(0,22){\circle*{4.1}}
\put(-15,22){\circle*{4.1}}
\put(15,22){\circle*{4.1}}
\put(30,22){\circle*{4.1}}
\put(34,25){$y$}
\put(0,47){\circle*{4.1}}
\put(-30,47){\circle*{4.1}}
\put(-15,47){\circle*{4.1}}
\put(15,47){\circle*{4.1}}
\put(30,47){\circle*{4.1}}
\put(30,47){\circle*{4.1}}
\put(-30,47){\line(1,0){59}}

\put(-30,-5){\circle*{4.1}}
\put(-42,-5){$y_1$}
\put(-30,-5){\line(1,0){59}}

\put(-44,49){$x_1$}
\put(35,49){$x_r$}

\put(0,-5){\circle*{4.1}}
\put(-15,-5){\circle*{4.1}}
\put(15,-5){\circle*{4.1}}
\put(30,-5){\circle*{4.1}}
\put(30,-5){\line(0,1){52}}
\put(35,-5){$y_s$}

\put(-7,35){$c_1$}

\put(-7,5){$c_2$}

\put(-7,53){$c_3$}
\end{picture}

\noindent Observe that, by the choice of the ${\cal P}_i$'s, a chord
of the cycle $c_1$ (resp. $c_2$, $c_3$) must join $x_i$ and $z_j$
(resp. $z_i$ and $y_j$, $x_i$ and $y_j$) for some $i,j$. 
Notice that the cycles $c_1$ and $c_3$ are primitive. Indeed
if $c_1$ or $c_3$ have a chord, then one of the following 
$$
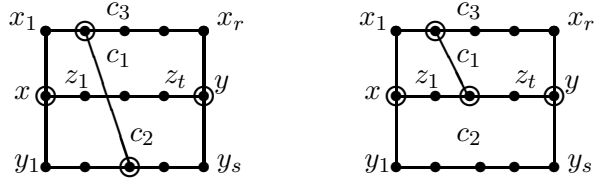
\begin{array}{cccc}
\begin{picture}(80,100)(-150,-25)
\setlength{\unitlength}{.035cm}
\thicklines
\put(-30,22){\circle*{4.1}}
\put(-30,22){\line(1,0){59}}
\put(-30,22){\circle{7.1}}
\put(-42,21){$x$}
\put(-30,47){\circle*{4.1}}
\put(-23,27){$z_1$}
\put(-30,21){\line(0,-1){26}}
\put(-30,21){\line(0,1){26}}
\put(15,27){$z_t$}
\put(0,22){\circle*{4.1}}
\put(-15,22){\circle*{4.1}}
\put(15,22){\circle*{4.1}}
\put(30,22){\circle*{4.1}}
\put(30,22){\circle{7.1}}
\put(34,25){$y$}
\put(0,47){\circle*{4.1}}
\put(-30,47){\circle*{4.1}}
\put(-15,47){\circle{7.1}}
\put(-15,47){\line(1,-3){17}}
\put(-15,47){\circle*{4.1}}
\put(15,47){\circle*{4.1}}
\put(30,47){\circle*{4.1}}
\put(30,47){\circle*{4.1}}
\put(-30,47){\line(1,0){59}}

\put(-30,-5){\circle*{4.1}}
\put(-42,-5){$y_1$}
\put(-30,-5){\line(1,0){59}}

\put(-44,49){$x_1$}
\put(35,49){$x_r$}
\put(2,-5){\circle{7.1}}
\put(2,-5){\circle*{4.1}}
\put(-15,-5){\circle*{4.1}}
\put(15,-5){\circle*{4.1}}
\put(30,-5){\circle*{4.1}}
\put(30,-5){\line(0,1){52}}
\put(35,-5){$y_s$}

\put(-7,35){$c_1$}

\put(2,5){$c_2$}

\put(-7,53){$c_3$}

\end{picture}
&\ \ \ \ \ \ \ \ \  &
\begin{picture}(80,100)(-150,-25)
\setlength{\unitlength}{.035cm}
\thicklines
\put(-30,22){\circle*{4.1}}
\put(-30,22){\line(1,0){59}}
\put(-30,22){\circle{7.1}}
\put(-42,21){$x$}
\put(-30,47){\circle*{4.1}}
\put(-23,27){$z_1$}
\put(-30,21){\line(0,-1){26}}
\put(-30,21){\line(0,1){26}}
\put(15,27){$z_t$}
\put(-2,22){\circle*{4.1}}
\put(-2,22){\circle{7.1}}
\put(-15,22){\circle*{4.1}}
\put(15,22){\circle*{4.1}}
\put(30,22){\circle*{4.1}}
\put(30,22){\circle{7.1}}
\put(34,25){$y$}
\put(0,47){\circle*{4.1}}
\put(-30,47){\circle*{4.1}}
\put(-15,47){\circle{7.1}}
\put(-15,47){\line(1,-2){13}}
\put(-15,47){\circle*{4.1}}
\put(15,47){\circle*{4.1}}
\put(30,47){\circle*{4.1}}
\put(30,47){\circle*{4.1}}
\put(-30,47){\line(1,0){59}}

\put(-30,-5){\circle*{4.1}}
\put(-42,-5){$y_1$}
\put(-30,-5){\line(1,0){59}}

\put(-44,49){$x_1$}
\put(35,49){$x_r$}
\put(2,-5){\circle*{4.1}}
\put(-15,-5){\circle*{4.1}}
\put(15,-5){\circle*{4.1}}
\put(30,-5){\circle*{4.1}}
\put(30,-5){\line(0,1){52}}
\put(35,-5){$y_s$}

\put(-7,35){$c_1$}

\put(-7,5){$c_2$}

\put(-7,53){$c_3$}

\end{picture}
&\ \ \ \ \ \ \ \ \ \ \ \ \ \ \ \ \ \ \ \ \ \ \ \ \ \ \ \ \ \ 
\ \ \ \ \ \ \ \ \ \ \ \  
\end{array}
$$
is a subgraph of $G$, which is impossible because both subgraphs 
are subdivisions of $\mathcal{K}_4$. Since $c_1$ and $c_3$ are
primitive and have at least two edges in common we obtain that $G$
does not satisfy PCP, a contradiction. \QED

\begin{Lemma}\label{addlemma1} Let $G$ be a graph. If ${\rm
rank}(G)={\rm frank}(G)$, then $G$ does not contain a subdivision 
of ${\cal K}_4$ as a subgraph. 
\end{Lemma}

\demo Assume there is a subgraph $H\subset G$
which is a subdivision of ${\cal K}_4$. If ${\cal K}_4$ is a subgraph
of $G$, then $G$ has four distinct triangles whose cycle vectors are
linearly dependent, a contradiction to Corollary~\ref{vila-isi}. If 
${\cal K}_4$ is not a subgraph of $G$, then $H$ is a strict subdivision
of ${\cal K}_4$, i.e., $H$ has more than four vertices. 
It follows that there are two vertices
$x,y$ in $V(H)$ which are non adjacent in $G$. Notice that $x,y$ 
can be chosen in ${\cal K}_4$ before subdivision. 
Therefore there are at least three non
adjacent paths joining $x$ and $y$, a contradiction to
Lemma~\ref{lina}. 
\QED

\medskip

The main result of this section is:

\begin{Theorem}\label{isidoro-vila05} Let $G$ be a graph. Then the 
following conditions are equivalent{\rm :}\vspace{-2mm}
\begin{description}
\item{\rm(a)} $G$ is a ring graph.\vspace{-2mm}
\item{\rm(b)} ${\rm rank}(G)={\rm frank}(G)$.\vspace{-2mm}
\item{\rm(c)} $G$ satisfies {\rm PCP} and $G$ does not
contain a subdivision of ${\cal K}_4$ as a subgraph. 
\end{description}
\end{Theorem}

\demo (a) $\Rightarrow$ (b): By induction on the number of vertices it is not hard 
to see that any ring graph $G$ satisfies the 
equality ${\rm rank}(G)={\rm frank}(G)$. 

(b) $\Rightarrow$ (c):  It follows at once from Lemmas~\ref{march11-03} and 
\ref{addlemma1}.

(c) $\Rightarrow$ (a): Let $G_1,\ldots,G_r$ be the blocks of $G$. 
The proof is by induction on the number of 
vertices of $G$. If each $G_i$ is
either a bridge or an isolated vertex, then $G$ is a forest and
consequently a ring graph. Hence by Lemma~\ref{mar22-05} we may
assume that $G$ is $2$-connected and that $G$ is not a cycle. We claim that $G$ 
has at least one vertex of degree $2$. If $\deg(v)\geq 3$ for all $v\in
V(G)$, then by Lemma~\ref{dirac05} there is a subgraph $H\subset G$
which is a subdivision of ${\cal K}_4$, which is impossible. 
Let $v_0\in V(G)$ be a vertex of degree 
$2$ as claimed. By the primitive cycle property 
there is a unique primitive cycle 
$c=\{v_0,v_1,\ldots,v_s=v_0\}$ of $G$ containing  $v_0$. The graph 
$H=G\setminus\{v_0\}$ satisfies PCP and does not has a subdivision 
of $\mathcal{K}_4$ as a subgraph. Consequently $H$ is a ring graph.
Thus we may assume that $c$ is not a triangle, otherwise $G$ is a ring
graph because it can be obtained by adding the $H$-path
$\{v_2,v_0,v_1\}$ to $H$. 

Next we claim that if $1\leq i<j<k\leq s-1$, then $v_i$ and
$v_k$ cannot be in the same connected component of
$H\setminus\{v_j\}$. Otherwise there is a path of $H\setminus\{v_j\}$
than joins $v_i$ with $v_k$. It follows that there is a path 
$\cal P$ of $H\setminus\{v_j\}$ with at least three vertices that joins a vertex of
$\{v_{j+1},\ldots,v_{s-1}\}$ with a vertex of
$\{v_{1},\ldots,v_{j-1}\}$ and such that $\cal P$ intersects $c$ 
exactly in its ends, but this contradicts Lemma~\ref{addlemma2}. This
proves the claim. In
particular $v_i$ is a cutvertex of $H$ for $i=2,\ldots,s-2$ and
$v_{i-1}$, $v_{i+1}$ are  
in different connected components of $H\setminus\{v_i\}$. For each
$1\leq i\leq s-2$ there is a block $K_i$ of $H$ such that
$\{v_i,v_{i+1}\}$ is an 
edge of $K_i$. Notice that if $1\leq i<j<k\leq s-1$, then 
$v_i,v_j,v_k$ cannot lie in some $K_\ell$. Indeed if the three vertices lie in
some $K_\ell$, then there is a path ${\cal P}'$ in $K_\ell\setminus\{v_j\}$ that
joins $v_i$ and $v_k$. Since ${\cal P}'$ is also a path in 
$H\setminus\{v_j\}$, we get that $v_i$ and $v_k$ are in the same
connected component of $H\setminus\{v_j\}$, but this contradicts the
last claim. 
In particular 
$V(K_\ell)$ intersects the cycle $c$ in exactly the vertices $v_\ell,v_{\ell+1}$ for
$1\leq \ell\leq s-2$.   


Observe that at least one of the edges of $c$ not containing $v_0$ is 
not a bridge of $H$. To show this pick $x\notin c$ such that $\{x,v_k\}$
is an edge of $H$. We may assume that $v_{k+1}\neq v_0$ 
(or $v_{k-1}\neq v_0$). Since $G'=G\setminus\{v_k\}$ is connected, there
is a path $\cal P$ of $G'$ joining $x$ and $v_{k+1}$ (or $v_{k-1}$).
This readily yields a  cycle of $H$ containing an edge of $c$ which
is not a bridge of 
$H$. Hence at least one of the blocks $K_1,\ldots,K_{s-2}$, 
say $K_i$, contains vertices outside $c$.

Next we show that two distinct blocks $B_1,B_2$ of $H$ cannot
intersect outside $c$. We 
proceed by contradiction assuming that $V(B_1)\cap V(B_2)=\{z\}$ for
some $z$ not in $c$. Let $H_1,\ldots,H_t$ be the connected components
of $H\setminus\{z\}$. Notice that $t\geq 2$ because $\{z\}$ is the
intersection of two different blocks of $H$. We may assume that
$\{v_1,\ldots,v_{s-1}\}$ are 
contained in $H_1$. Consider the subgraph $H_1'$ of $G\setminus\{z\}$
obtained from 
$H_1$ by adding the vertex $v_0$ and the edges $\{v_0,v_1\}$, 
$\{v_0,v_{s-1}\}$. It follows that the  connected
components of $G\setminus\{z\}$ are $H_1',H_2,\ldots,H_t$ , which is 
impossible because $G$ is
$2$-connected. 

Let $K_i$ be a block of $H$ that contains 
vertices outside $c$ for some $1\leq i\leq s-2$. By induction
hypothesis $K_i$ is a ring graph. Thus by 
Remark~\ref{march23-love} we can construct $K_i$ starting with a
primitive cycle $c_1$ that contains the edge $\{v_i,v_{i+1}\}$, and then
adding appropriate paths. Suppose that ${\cal P}_1,\ldots,{\cal P}_m$
is the sequence of paths added to $c_1$ to obtain $K_i$. If we remove the
path ${\cal P}_m$ from $G$ and use the fact that distinct blocks of
$H$ cannot intersect outside $c$, 
then again by induction hypothesis we obtain a ring
graph. It follows that $G$ is a ring graph as well. \QED

\medskip

An immediate consequence of Theorem~\ref{isidoro-vila05} is: 
\begin{Corollary}\label{march21-2-05} Let $G$ be a graph. If ${\rm rank}(G)={\rm
frank}(G)$, then $G$ is planar. 
\end{Corollary}

\begin{Corollary}\label{ring-induced} If $G$ is a ring graph and $H$
is an induced subgraph of $G$, then $H$ is a ring graph. 
\end{Corollary}

\demo It follows from part (c) of Theorem~\ref{isidoro-vila05}. \QED

\medskip



Two graphs $H_1$ and $H_2$ are called {\it
homeomorphic}
if there exists a graph $G$ such that both
$H_1$ and $H_2$ are subdivisions of $G$. A graph is {\it
outerplanar\/} if it can be embedded 
in the plane so that all its vertices lie on a common face; it is usual to
choose this face to be the exterior face. The complete bipartite 
graph with bipartition $(V_1,V_2)$ is denoted by ${\cal K}_{t,s}$,
where $|V_1|=t$ and $|V_2|=s$.

\begin{Theorem}{\rm \cite[Theorem~11.10]{Har}}\label{outerplanar-desc}
A graph is outerplanar if and only if it has no 
subgraph homeomorphic to ${\cal K}_4$ or ${\cal K}_{2,3}$ 
except 
${\cal K}_4\setminus\{e\}$, where $e$ is an edge.
\end{Theorem}

\begin{Proposition}\label{march11-1-03} If $G$ is an 
outerplanar graph, then ${\rm rank}(G)={\rm frank}(G)$.
\end{Proposition}

\demo By Theorem~\ref{isidoro-vila05}(c) it suffices to 
prove that $G$ satisfies {\rm PCP} and $G$ does not 
contain a subdivision of ${\cal K}_4$ as a subgraph. 
If $G$ contains a subdivision $H$ of ${\cal K}_4$ as a subgraph, 
then $G$ contains a subgraph, namely $H$, homeomorphic 
to ${\cal K}_4$, but this is impossible by
Theorem~\ref{outerplanar-desc}. To finish the proof we now show that 
$G$ has the PCP property. Let $c_1=\{x_1,x_2,\ldots,x_m=x_1\}$ and 
$c_2=\{y_1,y_2,\ldots,y_n=y_1\}$ be two distinct primitive cycles
having at least one common edge. We may assume that $x_i=y_i$ for 
$i=1,2$ and $x_3\neq y_3$. Notice that $y_3\notin c_1$ because
otherwise $\{y_2,y_3\}=\{x_2,y_3\}$ is a chord of $c_1$. We need only
show that $\{x_1,x_2\}=c_1\cap c_2$, because this implies that 
$c_1$ and $c_2$ cannot have more than one edge in common. 
Assume that
$\{x_1,x_2\}\subsetneq c_1\cap c_2$. Let $r$ be the minimum integer 
such that $y_r$ belong to $(c_1\cap c_2)\setminus\{x_1,x_2\}$. Notice
that $y_r\neq x_3$ because otherwise $\{x_2,x_3\}$ is a chord of
$c_2$. Hence $c_1$ together with the path 
$\{x_2=y_2,y_3,\ldots,y_r\}$ give a subgraph $H$ of $G$ which is a 
subdivision of ${\cal K}_{2,3}$, a contradiction to
Theorem~\ref{outerplanar-desc}. 
\QED

\section{Toric ideals of graphs}\label{simplegraphs}

Let $R=k[x_1,\ldots,x_n]$ be a polynomial ring 
over a field $k$ and let $G$ be a graph on the vertex 
set $V(G)=\{x_1,\ldots,x_n\}$.  The {\it edge subring\/} of the 
graph $G$, denoted by $k[G]$, is the $k$-subalgebra of $R$ generated 
by the monomials corresponding to the edges of $G$:
\[
k[G]=k[\{x_ix_j\vert\, x_i\, \mbox{ is adjacent to }x_j\}]
\subset R.
\]
There is a
graded 
epimorphism of $k$-algebras
\[ 
\varphi\colon B=k[t_1,\ldots,t_q] \longrightarrow k[G], \ \ \ 
\{x,y\}\longmapsto xy,
 \]
where $B$ is a polynomial ring graded by $\deg(t_i)=1$ for all 
$i$ and $k[G]$ has the normalized grading $\deg(f_i)=1$ for all $i$. 
The kernel of $\varphi$, denoted by $P(G)$, is a graded 
prime ideal 
of $B$ called the {\it toric ideal\/} 
of $G$. The graded structure of $P(G)$ will not play a role in what
follows. Later we will emphasize the fact that toric ideals of 
oriented graphs may not have a graded
structure. Having a grading is useful if one studies the projective 
toric variety defined by $P(G)$ or systems of generators of $P(G)$.   

The Krull dimension of $k[G]$ equals the rank of the incidence 
matrix of $G$ \cite{He3}. If $G$ is a connected 
graph, then by \cite[Corollary~6.3]{handbook} one has:
$$
\dim(k[G])=\left\{\begin{array}{ll}
n&\mbox{if }G \mbox{ is not bipartite}, \mbox{ and}\\ 
n-1&\mbox{otherwise}.
\end{array} \right. 
$$
Since $B/P(G)\simeq k[G]$, we obtain
that height of $P(G)$ is $q-n+1$ 
if $G$ is  
a connected bipartite graph and that height of $P(G)$ is $q-n$ if $G$ 
is a connected non-bipartite graph.

\begin{Definition}\rm The toric ideal $P(G)$ is
called a {\it complete
intersection\/} if it can be generated by $g$ polynomials, where $g$ is the 
height of $P(G)$. The graph $G$ is called a {\it complete 
intersection\/} if $P(G)$ is a complete
intersection.  
\end{Definition}

The complete intersection property is independent of $k$
\cite[Theorem~3.9]{morales-thoma}. In
the area of complete intersection toric ideals there are some 
recent papers, see \cite{stcib,stcib-algorithm} 
and the introduction of \cite{morales-thoma}, where one can find 
additional properties and references on this active area. 

Next we describe a generating set for $P(G)$ that shows how the cycle
structure of $G$ determine $P(G)$. Let 
$$c=\{x_0,x_1,\ldots,x_r=x_0\}$$ 
be an even cycle of $G$ such 
that $f_i=x_{i-1}x_i$. Notice that the binomial 
$$t_c=t_1t_3\cdots
t_{r-1}-t_2t_4\cdots t_r$$ 
is in $P(G)$. If $G$ is 
bipartite, then $P(G)$ is minimally generated by the set of all $t_c$
such that $c$ is a primitive cycle of $G$, see \cite{Vi3}.  

\medskip

The next result can be 
extended to non connected
bipartite graphs. 

\begin{Theorem}{\rm\cite[Theorem~2.5]{aron-jac}}
\label{aronjacp} If $G$ is a
bipartite connected graph, then $G$ is a 
complete intersection if and only if ${\rm rank}(G)={\rm frank}(G)$.
\end{Theorem}

This was the first characterization of complete intersection bipartite
graphs. For these graphs the equality ${\rm rank}(G)={\rm frank}(G)$ can also be 
interpreted in homological terms \cite{aron-jac}.
Another characterization is the following:

\begin{Theorem}[\rm\cite{katzman}]\label{mar15-01} If $G$ is a 
bipartite graph, then $G$ is a complete intersection 
if and only if $G$ is planar and satisfies {\rm PCP}. 
\end{Theorem}

The next result 
is interesting because it shows how to construct all the complete intersection
bipartite graphs. 

\begin{Corollary}\label{aug20-05} If $G$ is a bipartite graph, then $G$ is a complete
intersection if and only if $G$ is a ring graph. 
\end{Corollary}

\demo By Theorem~\ref{aronjacp} $G$ is a complete intersection 
if and only if ${\rm rank}(G)={\rm frank}(G)$ and the result follows 
from Theorem~\ref{isidoro-vila05}. \QED

\medskip

\noindent{\it Notation\/} For 
$a=({a}_1,\ldots,{a}_q)\in {\mathbb N}^q$ and
$f_1,\ldots, f_q$ in a commutative ring we 
set $f^a
= f_1^{a_1}\cdots f_q^{a_q}$. The {\it support\/}
 of $f^a$ is the set ${\rm supp}(f^a)= \{f_i\, |\,
a_i\neq 0\}$. 

\begin{Definition}\rm Let
$g_1=t^{\alpha_1}-t^{\beta_1},\ldots,g_r=t^{\alpha_r}-t^{\beta_r}$ be
a sequence of homogeneous binomials of degree at least $2$ in the polynomial ring
$B=k[t_1,\ldots,t_q]$. We say  that ${\cal B}=\{g_1,\ldots,g_r\}$ is
a {\it foliation\/} if the following conditions are
satisfied:\vspace{-1mm}
\begin{description}
\item{\rm (a)} $t^{\alpha_i}$ and $t^{\beta_i}$ are square-free 
monomials for
all $i$,\vspace{-2mm}
\item{\rm (b)} ${\rm supp}(t^{\alpha_i})\cap{\rm
supp}(t^{\beta_i})=\emptyset$ for all $i$, and\vspace{-2mm}
\item{(c)} $|(\cup_{i=1}^j C_i)\cap C_{j+1}|=1$ for  $1\leq j<r$,
where $C_i={\rm supp}(t^{\alpha_i})\cup{\rm
supp}(t^{\beta_i})$. 
\end{description}
\end{Definition}

\begin{Proposition}\label{foliation} If ${\cal B}=\{g_1,\ldots,g_r\}$
is a foliation, 
then the ideal $I=({\cal B})$ generated by $\mathcal{B}$ 
is a complete intersection and
${\cal B}$ is a Gr\"obner basis of $I$.
\end{Proposition}

\demo By the constructive nature of ${\cal B}$ we can order the
variables $t_1,\ldots,t_q$ such that the leading terms of
$g_1,\ldots,g_r$, with respect to the lexicographical order, 
are relatively prime. Let ${\rm in}(g_i)$ be the leading term of
$g_i$. Then ${\cal B}$ is
a Gr\"obner basis by a result of 
Buchberger \cite[Theorem~2.4.15]{monalg}. Since $B/I$ and $B/({\rm
in}(g_1),\ldots,{\rm in}(g_r))$  
have the same Krull dimension by a result of Macaulay
\cite[Corollary~2.4.13]{monalg}, we obtain that the height of $I$ is
equal to $r$, as required. \QED 
 
\begin{Corollary}\label{jul16-05} If $G$ is a $2$-connected bipartite graph with at
least four vertices, then the toric ideal $P(G)$ is a complete
intersection if and only if it is generated by a foliation. 
\end{Corollary}

\demo It follows from Corollary~\ref{aug20-05} and the definition of a
 ring graph. \QED

\section{Toric ideals of oriented graphs}\label{orientedgraphs}

Let $G$ be a connected graph with $n$ vertices and $q$ edges and let $\cal O$ be an 
orientation of the edges of $G$, i.e., an 
assignment of a direction to each edge of $G$. 
Thus ${\cal D}=(G,{\cal O})$ is an {\it oriented
graph}. To each oriented edge $e=(x_i,x_j)$ of $\cal D$, we associate
the vector $v_e$ defined as follows: the $i${\it th} entry 
is $-1$, the $j${\it th} entry is $1$, and the remaining entries are zero. The
{\it incidence matrix\/} $A_{\cal D}$ of $\cal D$ is
the $n\times q$ matrix with entries in $\{0,\pm 1\}$ whose 
columns are the vectors of the form $v_e$, with $e$ an edge of $\cal
D$. For simplicity of notation we set $A=A_{\cal D}$. 
The set of column vectors of $A$ will be denoted by 
${\cal A}=\{v_1,\ldots,v_q\}$. It is well known \cite{oxley}
that $A$ defines a matroid $M[A]$ on ${\cal A}=\{v_1,\ldots,v_q\}$ over the field
$\mathbb{Q}$ of rational numbers, which is called the {\it vector
matroid\/} of $A$, whose independent sets are
the independent subsets of ${\cal A}$. A {\it minimal dependent set\/} or {\it
circuit\/} of $M[A]$ is a
dependent set all of whose proper subsets are independent. A subset
$B$ of ${\cal A}$ is called a {\it basis\/} of
$M[A]$ if $B$ is a maximal
independent set. Recall that an integer matrix is called {\it totally
unimodular\/} if each $i\times i$ minor (subdeterminant) of the matrix
 is $0$ or $\pm 1$ for all $i\geq 1$.

\begin{Lemma}\label{jul17-05}
The circuits of $M[A]$ are precisely 
the cycles of $G$, $A$ is totally unimodular, and 
${\rm rank}(A)=n-1$
\end{Lemma}

\demo It follows from \cite[pp.~343-344]{godsil} and \cite[p.
274]{Schr}. \QED

\medskip

Let $\alpha\in {\mathbb R}^q$. The {\it support\/} of $\alpha$ 
is defined as ${\rm supp}(\alpha)=\{i\, |\, \alpha_i\neq
0\}$. An {\it elementary
vector\/} of ${\rm ker}(A)$ is a vector $0\neq \alpha$ in ${\rm
ker}(A)$ whose support is minimal with 
respect to inclusion, i.e., ${\rm supp}(\alpha)$ does not properly
contain the support of any other nonzero vector in ${\rm ker}(A)$. A 
{\it circuit\/} of ${\rm ker}(A)$ 
is an elementary vector of ${\rm ker}(A)$ 
with relatively prime integral entries (see \cite[Section~2]{accota}).  
There is a one to one correspondence
$$
\begin{array}{ccc}
\mbox{Circuits of } \ker(A)&\longrightarrow & \mbox{ Circuits of }
M[A]=\mbox{cycles of }G
\end{array}
$$
given by $\alpha=(\alpha_1,\ldots,\alpha_q)\rightarrow 
C(\alpha)=\{v_i\vert\, i\in{\rm supp}(\alpha)\}$. 
Thus the set of circuits of the kernel of $A$ is the
algebraic realization of the set of circuits of the vector matroid
$M[A]$.

Consider the {\it edge subring\/} $k[{\cal D}]:=k[x^{v_1},\ldots,x^{v_q}]
\subset k[x_1^{\pm 1},\ldots,x_n^{\pm 1}]$ of the oriented graph 
${\cal D}$. There is an epimorphism of $k$-algebras
\[ 
\varphi\colon B=k[t_1,\ldots,t_q] \longrightarrow k[{\cal D}], \ \ \ 
t_i\longmapsto x^{v_i},
 \]
where $B$ is a polynomial ring. The kernel of $\varphi$, denoted by
$P_{\cal D}$, is called the {\it toric ideal\/} of 
${\cal D}$. Notice that $P_{\cal
D}$ is no longer a graded ideal, see Proposition~\ref{feb21-06}.
The toric ideal $P_{\cal D}$ is a prime ideal of height $q-n+1$
generated by binomials and $k[{\cal D}]$ is a normal domain. Thus
any minimal generating set of $P_{\cal D}$ must have at least $q-n+1$
elements, by the principal ideal theorem.

Let $\alpha\in\mathbb{R}^q$. Note that $\alpha=\alpha_+-\alpha_-$, 
where $\alpha_+$ and 
$\alpha_-$ are two non negative vectors 
with disjoint support. If $0\neq\alpha\in{\rm ker}(A)\cap\mathbb{Z}^n$ we associate the 
binomial $t_\alpha=t^{\alpha_+}-t^{\alpha_-}$. 
Notice that $t_\alpha\in P_{\cal D}$. Given a cycle $c$ of ${\cal
D}$, we split $c$ in two disjoint sets of edges $c_+$ 
and $c_-$, where $c_+$ is oriented clockwise and $c_-=c\setminus c_+$. The binomial 
$$
t_c=\prod_{v_i\in c_+}t_i-\prod_{v_i\in c_-}t_i
$$
belongs to $P_{\cal D}$. If $c_+=\emptyset$ or $c_-=\emptyset$ we set
$\prod_{v_i\in c_+}t_i=1\ \ \mbox{or}\ \ \prod_{v_i\in c_-}t_i=1.$

\begin{Definition}\rm The toric ideal $P_{\cal D}$ is called a 
{\it binomial complete intersection\/} if $P_{\cal D}$ can be generated by $q-n+1$
binomials. 
\end{Definition}

If $P_{\cal D}$ is homogeneous 
and is generated by $q-n+1$ polynomials, then 
$P_{\cal D}$ is a binomial complete intersection.

\begin{Proposition}\label{grobnbas} 
$P_{\cal D}$ is generated by the set of all binomials
$t_c$ such that $c$ is a cycle of ${\cal D}$ and this set is a universal 
Gr\"obner basis.
\end{Proposition}

\demo  Let ${\cal U}_{\cal D}$ be the set of all binomials of the form
$t_\alpha$ such that $\alpha$ is a circuit of ${\rm ker}(A)$. 
Since $A$ is totally unimodular, by \cite[Proposition~8.11]{Stur1}, 
the set ${\cal U}_{\cal D}$ form a universal Gr\"obner basis of
$P_{\cal D}$. Notice that the circuits of ${\rm ker}(A)$ are in one to
one correspondence with the circuits of the vector matroid $M[A]$.
To complete the proof it suffices to observe that the circuits of
$M[A]$ 
are precisely the cycles of $G$, see Lemma~\ref{jul17-05}. \QED

\begin{Proposition}\label{cyclepri} 
Let $c=\{x_1,x_2,\ldots,x_r,x_1\}$ be a circuit of ${\cal D}$. Suppose
that $(x_i,x_j)$ or $(x_j,x_i)$ is an edge of ${\cal D}$, 
with $i + 1 < j$. Then $t_c$ is a linear
combination of $t_{c_1}$ and $t_{c_2}$, where 
$c_1=\{x_1,x_2,\dots,x_i,x_j,x_{j+1},\dots,x_r,x_1\}$
and $c_2=\{x_i,x_{i+1},\dots,x_j,x_i\}$.
\end{Proposition}

\demo 
Suppose without loss of generality that $v_k=(x_i,x_j)$ is the edge
of ${\cal D}$ with $i+1<j$. Then we can write 
$t_{c_1}=t^{\alpha_+}-t^{\alpha_-}$ and
$t_{c_2}=t^{\beta_+}-t^{\beta_-}$ for some $\alpha$, $\beta$. 
We may assume that $v_k \in c_{1_+}\cap c_{2_+}$, because 
otherwise we may multiply $t_{c_1}$ or $t_{c_2}$ by
$-1$. As $t_k$ divides $t^{\alpha_+}$ and $t_k$ divides
$t^{\beta_+}$, we get 
\begin{eqnarray*}
\left(\frac{t^{\beta_+}}{t_k}\right)t_{c_1}-\left(\frac{t^{\alpha_+}}{t_k}\right)
t_{c_2}&=&\left(\frac{t^{\beta_+}}{t_k}\right)
(t^{\alpha_+}-t^{\alpha_-})-\left(\frac{t^{\alpha_+}}{t_k}\right)
(t^{\beta_+}-t^{\beta_-})\\ 
\\
 &=&\left(\frac{t^{\alpha_+}}{t_k}\right)t^{\beta_-}-
\left(\frac{t^{\beta_+}}{t_k}\right) 
 t^{\alpha_-} = t^{\gamma_1} - t^{\gamma_2}.
\end{eqnarray*}
Hence $t^{\gamma_1} - t^{\gamma_2}$ is 
in $P_{\cal D}$, where $\gamma_1=({\alpha_+}-e_k)+\beta_-$ 
and $\gamma_2=({\beta_+}-e_k)+\alpha_-$. Then $t^{\gamma_1}$
is the product of the edges of $(c_{1_+}\setminus\{t_k\})\cup c_{2_-}$,
but these are the edges of $c_+$. By the same reason 
$t_{\gamma_2}$ is the product of the edges of $c_-$. Thus
$t_c=t^{\gamma_1}-t^{\gamma_2}$. From the equality above we get
that $t_c$ is a linear combination of $t_{c_1}$ and $t_{c_2}$. 
\QED

\medskip

As an immediate consequence of Propositions~\ref{grobnbas} and ~\ref{cyclepri} 
we get:
\begin{Corollary}\label{jul22-1-05} $P_{\cal D}$ is generated by the set of binomials
corresponding to primitive cycles.
\end{Corollary}

We say that a cycle $c$ of ${\cal D}$ is {\it oriented\/} if all the arrows
of $c$ are oriented in the same direction. If ${\cal D}$
does not have oriented cycles, we say that ${\cal D}$ 
is acyclic. 

\begin{Proposition}[\rm\cite{Har}]
${\cal D}$ is acyclic if and only if there is a linear ordering of the
vertices such that every edge of 
${\cal D}$ has the form $(x_i,x_j)$ with $i<j$.  
\end{Proposition}

The ordering of the last proposition is called a {\it topological
ordering\/}. 
The next result is not hard to prove.

\begin{Proposition}\label{feb21-06}
If $\cal D$ has a topological ordering, then $P_{\cal D}$ is generated
by homogeneous binomials with respect to the grading induced by 
$degree(t_k)=j-i$, where $t_k$ maps to $x_i^{-1}x_j$ and $(x_i,x_j)$
is an edge.
\end{Proposition}

\begin{Corollary}\label{feb23-06}
If ${\cal D}$ is acyclic, then $P_{\cal D}$
is a complete intersection if and only if $P_{\cal D}$ is generated
by $q-n+1$ binomials corresponding to primitive cycles.  
\end{Corollary}

\demo Since $P_{\cal D}$ is a graded ideal, it suffices to 
recall that all the homogeneous minimal sets 
of generators of $P_{\cal D}$ have the same number of elements. \QED

\medskip

In general the binomial complete intersection property of
$P_{\cal D}$ depends on the orientation of $G$. However we have:

\begin{Corollary}\label{aug10-05} If $G$ is a ring graph, then
$P_{\cal D}$ 
is a complete intersection for any orientation of $G$. The converse
holds if $G$ is bipartite. 
\end{Corollary}

\demo By Corollary~\ref{jul22-1-05}, $P_{\cal D}$ is generated by $q-n+1$ 
binomials.  To show the converse assume that $G$ is bipartite. Let
$(V_1,V_2)$ be a bipartition of $G$. Consider the  
oriented graph $\cal D$ obtained from $G$ by orienting all the edges 
of $G$ from $V_1$ to $V_2$, i.e., all the arrows of $G$ have 
tail at $V_1$ and head at $V_2$. Since every vertex of $\cal D$ is
either a source or a sink it follows that $P(G)=P_{\cal D}$. Hence
$P(G)$ is a complete intersection and $G$ is a ring graph by
Corollary~\ref{aug20-05}. \QED

\medskip

An interesting problem that remains unsolved is to characterize the
graphs with the property that  $P_{\cal D}$  is a binomial complete
intersection for all orientations of $G$. Apart from ring graphs 
it has been shown that complete graphs have this property
\cite{tesis-enrique,morf}.

\subsection*{A special orientation}

Let $G$ be a connected graph. Here we show that there is always
an orientation of $G$
such that $P_{\cal D}$ is a complete intersection generated by the
binomials that correspond to a cycle basis of a certain spanning tree of $G$.  

\begin{Definition}\rm 
Let $S$ be a set of vertices of a graph $G$. The 
{\it neighbor set\/} of $S$, denoted by 
$N_G(S)$ or simply by $N(S)$ if $G$ is understood, is the 
set of vertices of $G$ that are adjacent with at least one 
vertex of $S$.
\end{Definition}

\begin{Lemma}\label{feb20-06} 
If $H$ is a subgraph of a connected graph $G$ and 
$N_{G}(V(H)) \subset V(H)$, then $V(G)=V(H)$.
\end{Lemma}
\demo Fix a vertex $x\in V(H)$. Let $y\in V(G)$. Since $G$ is
connected, there is a path  ${\cal P}=\{b_1 = x,b_2,\ldots,b_\ell=y\}$
from $x$ to $y$. Using that 
$\{b_j,b_{j+1}\} \in E(G)$ for $1\leq j<\ell-1$ and that $b_1\in V(H)$, 
by induction we get 
that $b_j\in V(H)$ for all $j$. Thus $y\in V(H)$. \QED

\medskip

We begin by constructing a proper nested sequence 
$A_1,\ldots, A_m$ of
subtrees of $G$ labeled by $V(A_j)=\{y_1^j,\ldots,y_{r_j}^j\}$
such that $A_m$ is a spanning tree of $G$ and $V(A_i)\subsetneq
V(A_{i+1})$ for $i<m$. First we construct the sequence
$A_1,\ldots,A_m$ and then we show that it has the required
properties.  
Let $A_1$ be a path of $G$ maximal with respect to inclusion.
Set $V(A_1)=\{y_1^1,y_2^1,\ldots,y_{r_1}^1\}$. We define
$$
i_1=\max\{u \in \mathbb {N}\vert N_{G}(y_1^1,\ldots,y_u^1)\subset V(A_1)\},
$$
where $N_{G}(B)$ is the neighbor set of $B$. 
If $i_1=r_1$, then $N_G(V(A_1))\subset V(A_1)$ and by
Lemma~\ref{feb20-06} we get $V(A_1)=V(G)$, in this case $A_1$ is the 
required spanning tree and we set $m=1$. 
If $i_1<r_1$, we define $a_1 = y_{i_1+1}^1$. By induction we define 
the sequence of subgraphs $A_1,\ldots,A_m$. 
Suppose that $A_j$ has been defined, 
where $V(A_j)=\{y_1^j,\ldots,y_{r_j}^j\}$. We define
$$
i_j = \max\{u \in \mathbb{N} \vert N_{G}(y_1^j,\ldots,y_{u}^j) \subset V(A_j)\}.
$$
If
$i_j=r_j$, then by Lemma~\ref{feb20-06} we get $V(A_j)=V(G)$, in this
case we set $m=j$ and $A_1,\ldots,A_j$ is the desired sequence. 
If $i_j<r_j=\vert V(A_j)\vert$, we define $a_j=y_{i_j+1}^j$.  Let
${\cal L}_j$ 
be  a maximal path with respect to inclusion such that 
$V({\cal L}_j) \cap V(A_j)=\{a_j\}$ and 
$V({\cal L}_j)=\{z_1^j,z_2^j,\ldots,z_{s_j}^j=a_j\}$,
the final vertex of ${\cal L}_j$ is $a_j$. We define $A_{j+1}$ as follows: 
$V(A_{j+1})=V(A_j)\cup V({\cal
L}_j)=\{y_1^{j+1},\ldots,y_{r_j+s_j-1}^{j+1}\}$, where
\begin{equation}
y_i^{j+1}=\left\{\begin{array}{ll}
                  y_i^j & \mbox{if $i\leq i_j$},\\
                  z_{i-i_j}^j & \mbox{if $i_j+1\leq i\leq i_j+s_j$},\\
                  y_{i-s_j+1}^j & \mbox{if $i_j+s_j+1 \leq i \leq
		  r_j+s_j-1$},
                 \end{array}
          \right.    
\end{equation}
$E(A_{j+1})=E(A_j)\cup E({\cal L}_j)$, and $r_{j+1}=r_j+s_j-1$. 

\begin{Lemma}\label{obs1}
$i_{k+1} > i_k$ for $1\leq k \leq m-1$.
\end{Lemma}

\demo By construction $y_i^{k+1}=y_i^k$ for $1\leq i \leq i_k$ and 
$y_{i_k+1}^{k+1}=z_1^k$
(see Eq.(1)). By the maximality of ${\cal L}_j$ we have 
$$
N_G(y_1^{k+1},y_2^{k+1},\ldots,y_{i_k}^{k+1},y_{i_k+1}^{k+1})
\subset V(A_{k+1}),
$$
thus $i_{k+1}>i_k$ by definition of $i_{k+1}$. \QED

\medskip

Suppose that the process finish at step $m$, i.e., $i_m=r_m$. We now
prove that $A_1,\ldots,A_m$ has the required properties:
\begin{Lemma}\label{spantree}
$A_i$ is a tree for $1\leq i \leq m$ and $A_m$ is a spanning tree 
of $G$.
\end{Lemma}

\demo By induction on $i$. For $i=1$ the assertion is clear. Suppose that $A_i$ is a tree. 
Recall that ${\cal L}_i$ is a tree and $V({\cal L}_i)\cap V(A_i)=\{a_i\}$. 
On the other hand $V(A_{i+1})=V(A_i)\cup V({\cal L}_i)$ 
and $E(A_{i+1})=E(A_i)\cup E({\cal L}_i)$, then $A_{i+1}$ is connected and
does not has cycles. By Lemma~\ref{feb20-06} we get 
that $V(A_m)=V(G)$ and $A_m$ is a 
spanning tree. 
\QED

\paragraph*{Orientation of the tree $A_m$ and the graph $G$.}

Let $\tau=(A_m,{\cal O})$ be the oriented tree obtained 
from $A_m$ using the following orientation:
$$
(y_i^m ,y_j^m) \in E(\tau) \mbox{ if and only if } \{y_i^m,y_j^m\}
\in E(A_m) \mbox{ and } j > i. 
$$

By Lemma~\ref{spantree} we have $V(G)=V(A_m)=\{y_1^m ,y_2^m,\ldots,y_{r_m}^m\}$
and we orient $G$ to obtain the oriented graph 
${\cal D}=(G,{\cal O})$ in the following way:
$$
(y_i^m,y_j^m) \in {\cal D} \mbox{ if and only if } \{y_i^m,y_j^m\} \in E(G) \mbox{ and } j>i.
$$

\begin{Example} 
The construction of the spanning tree $A_m$ and the orientation 
${\cal O}$ of $G$ is illustrated below. 

$$
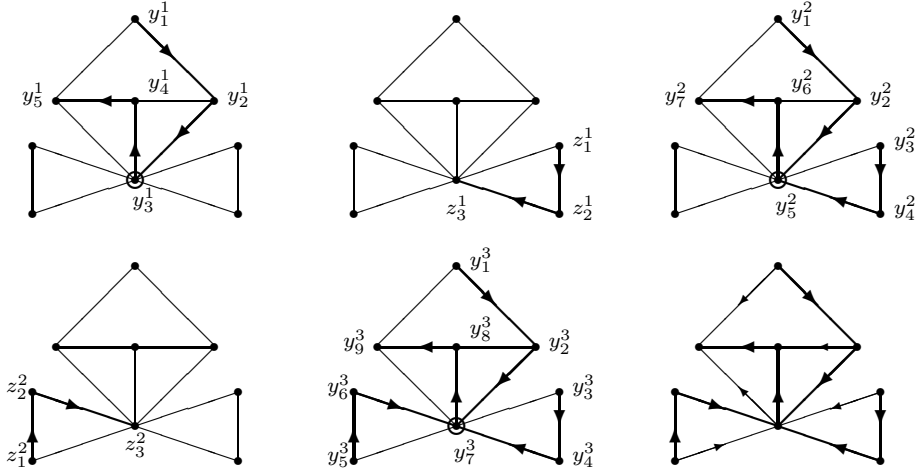
\begin{array}{cccccc}

\,\,\,\,\,\,\,\,\,\,\,\,\,\,\,\,\,\,\,\,\,\,\,
&
\thicklines
\begin{picture}(80,35)
\put(0,35){\circle*{3.1}}
\put(5,35){$\scriptstyle y_1^1$}

\put(0,35){\thicklines\line(1,-1){30}}
\put(0,35){\vector(1,-1){15}}
\put(30,4){\circle*{3.1}}
\put(35,4){$\scriptstyle y_2^1$}
\put(30,4){\thicklines\line(-1,-1){30}}
\put(30,4){\vector(-1,-1){15}}

\put(0,35){\thinlines\line(-1,-1){30}}
\put(-30,4){\circle*{3.1}}
\put(-43,4){$\scriptstyle y_5^1$}
\put(-30,4){\thinlines\line(1,-1){30}}

\put(0,-26){\circle*{3.1}}
\put(0,-26){\circle{6}}
\put(-1,-36){$\scriptstyle y_3^1$}

\put(0,-26){\thicklines\line(0,1){30}}
\put(0,-26){\vector(0,1){15}}
\put(0,4){\circle*{3.1}}
\put(5,9){$\scriptstyle y_4^1$}
\put(0,4){\thinlines\line(1,0){30}}
\put(0,4){\thicklines\line(-1,0){30}}
\put(0,4){\vector(-1,0){15}}

\put(0,-26){\thinlines\line(3,1){38}}
\put(0,-26){\thinlines\line(3,-1){38}}
\put(39,-13){\circle*{3.1}}
\put(39,-39){\circle*{3.1}}
\put(39,-13){\thinlines\line(0,-1){26}}

\put(0,-26){\thinlines\line(-3,1){38}}
\put(0,-26){\thinlines\line(-3,-1){38}}
\put(-39,-13){\circle*{3.1}}
\put(-39,-39){\circle*{3.1}}
\put(-39,-13){\thinlines\line(0,-1){26}}
\end{picture}

&
\,\,\,\,\,\,\,\,\,\,\,\,
&

\begin{picture}(80,35)
\put(0,35){\circle*{3.1}}

\put(0,35){\thinlines\line(1,-1){30}}
\put(30,4){\circle*{3.1}}
\put(30,4){\thinlines\line(-1,-1){30}}
 
\put(0,35){\thinlines\line(-1,-1){30}}
\put(-30,4){\circle*{3.1}}
\put(-30,4){\thinlines\line(1,-1){30}}

\put(0,-26){\circle*{3.1}}
\put(-4,-39){$\scriptstyle z_3^1$}

\put(0,-26){\thinlines\line(0,1){30}}
\put(0,4){\circle*{3.1}}
\put(0,4){\thinlines\line(1,0){30}}
\put(0,4){\thinlines\line(-1,0){30}}

\put(0,-26){\thinlines\line(3,1){38}}
\put(0,-26){\thicklines\line(3,-1){38}}
\put(39,-39){\vector(-3,1){19}}
\put(39,-13){\circle*{3.1}}
\put(44,-13){$\scriptstyle z_1^1$}
\put(39,-39){\circle*{3.1}}
\put(44,-39){$\scriptstyle z_2^1$}
\put(39,-13){\thicklines\line(0,-1){26}}
\put(39,-13){\vector(0,-1){13}}

\put(0,-26){\thinlines\line(-3,1){38}}
\put(0,-26){\thinlines\line(-3,-1){38}}
\put(-39,-13){\circle*{3.1}}
\put(-39,-39){\circle*{3.1}}
\put(-39,-13){\thinlines\line(0,-1){26}}
\end{picture}

&
\,\,\,\,\,\,\,\,\,\,\,\,
&

\begin{picture}(80,35)
\put(0,35){\circle*{3.1}}
\put(5,35){$\scriptstyle y_1^2$}

\put(0,35){\thicklines\line(1,-1){30}}
\put(0,35){\vector(1,-1){15}}
\put(30,4){\circle*{3.1}}
\put(35,4){$\scriptstyle y_2^2$}
\put(30,4){\thicklines\line(-1,-1){30}}
\put(30,4){\vector(-1,-1){15}}

\put(0,35){\thinlines\line(-1,-1){30}}
\put(-30,4){\circle*{3.1}}
\put(-43,4){$\scriptstyle y_7^2$}
\put(-30,4){\thinlines\line(1,-1){30}}

\put(0,-26){\circle*{3.1}}
\put(0,-26){\circle{6}}
\put(-1,-38){$\scriptstyle y_5^2$}

\put(0,-26){\thicklines\line(0,1){30}}
\put(0,-26){\vector(0,1){15}}
\put(0,4){\circle*{3.1}}
\put(5,9){$\scriptstyle y_6^2$}
\put(0,4){\thinlines\line(1,0){30}}
\put(0,4){\thicklines\line(-1,0){30}}
\put(0,4){\vector(-1,0){15}}

\put(0,-26){\thinlines\line(3,1){38}}
\put(0,-26){\thicklines\line(3,-1){38}}
\put(39,-39){\vector(-3,1){19}}
\put(39,-13){\circle*{3.1}}
\put(44,-13){$\scriptstyle y_3^2$}
\put(39,-39){\circle*{3.1}}
\put(44,-39){$\scriptstyle y_4^2$}
\put(39,-13){\thicklines\line(0,-1){26}}
\put(39,-13){\vector(0,-1){13}}

\put(0,-26){\thinlines\line(-3,1){38}}
\put(0,-26){\thinlines\line(-3,-1){38}}
\put(-39,-13){\circle*{3.1}}
\put(-39,-39){\circle*{3.1}}
\put(-39,-13){\thinlines\line(0,-1){26}}
\end{picture}

\\

& & & & &

\\

& & & & &

\\

& & & & &

\\

& & & & &

\\

\,\,\,\,\,\,\,\,\,\,\,\,\,\,\,\,\,\,\,\,

&

\begin{picture}(80,35)
\put(0,35){\circle*{3.1}}

\put(0,35){\thinlines\line(1,-1){30}}
\put(30,4){\circle*{3.1}}
\put(30,4){\thinlines\line(-1,-1){30}}

\put(0,35){\thinlines\line(-1,-1){30}}
\put(-30,4){\circle*{3.1}}
\put(-30,4){\thinlines\line(1,-1){30}}

\put(0,-26){\circle*{3.1}}
\put(-4,-35){$\scriptstyle z_3^2$}

\put(0,-26){\thinlines\line(0,1){30}}
\put(0,4){\circle*{3.1}}
\put(0,4){\thinlines\line(1,0){30}}
\put(0,4){\thinlines\line(-1,0){30}}

\put(0,-26){\thinlines\line(3,1){38}}
\put(0,-26){\thinlines\line(3,-1){38}}
\put(39,-13){\circle*{3.1}}
\put(39,-39){\circle*{3.1}}
\put(39,-13){\thinlines\line(0,-1){26}}

\put(0,-26){\thicklines\line(-3,1){38}}
\put(-39,-13){\vector(3,-1){19}}
\put(0,-26){\thinlines\line(-3,-1){38}}
\put(-39,-13){\circle*{3.1}}
\put(-49,-13){$\scriptstyle z_2^2$}
\put(-39,-39){\circle*{3.1}}
\put(-49,-39){$\scriptstyle z_1^2$}
\put(-39,-13){\thicklines\line(0,-1){26}}
\put(-39,-39){\vector(0,1){13}}
\end{picture}

&
\,\,\,\,\,\,\,\,\,\,\,\,
&

\begin{picture}(80,35)
\put(0,35){\circle*{3.1}}
\put(5,35){$\scriptstyle y_1^3$}

\put(0,35){\thicklines\line(1,-1){30}}
\put(0,35){\vector(1,-1){15}}
\put(30,4){\circle*{3.1}}
\put(35,4){$\scriptstyle y_2^3$}
\put(30,4){\thicklines\line(-1,-1){30}}
\put(30,4){\vector(-1,-1){15}}

\put(0,35){\thinlines\line(-1,-1){30}}
\put(-30,4){\circle*{3.1}}
\put(-43,4){$\scriptstyle y_9^3$}
\put(-30,4){\thinlines\line(1,-1){30}}

\put(0,-26){\circle*{3.1}}
\put(0,-26){\circle{6}}
\put(-1,-38){$\scriptstyle y_7^3$}

\put(0,-26){\thicklines\line(0,1){30}}
\put(0,-26){\vector(0,1){15}}
\put(0,4){\circle*{3.1}}
\put(5,9){$\scriptstyle y_8^3$}
\put(0,4){\thinlines\line(1,0){30}}
\put(0,4){\thicklines\line(-1,0){30}}
\put(0,4){\vector(-1,0){15}}

\put(0,-26){\thinlines\line(3,1){38}}
\put(0,-26){\thicklines\line(3,-1){38}}
\put(39,-39){\vector(-3,1){19}}
\put(39,-13){\circle*{3.1}}
\put(44,-13){$\scriptstyle y_3^3$}
\put(39,-39){\circle*{3.1}}
\put(44,-39){$\scriptstyle y_4^3$}
\put(39,-13){\thicklines\line(0,-1){26}}
\put(39,-13){\vector(0,-1){13}}

\put(0,-26){\thicklines\line(-3,1){38}}
\put(-39,-13){\vector(3,-1){19}}
\put(0,-26){\thinlines\line(-3,-1){38}}
\put(-39,-13){\circle*{3.1}}
\put(-49,-13){$\scriptstyle y_6^3$}
\put(-39,-39){\circle*{3.1}}
\put(-49,-39){$\scriptstyle y_5^3$}
\put(-39,-13){\thicklines\line(0,-1){26}}
\put(-39,-39){\vector(0,1){13}}
\end{picture}

&
\,\,\,\,\,\,\,\,\,\,\,\,
&

\begin{picture}(80,35)
\put(0,35){\circle*{3.1}}

\put(0,35){\thicklines\line(1,-1){30}}
\put(0,35){\vector(1,-1){15}}
\put(30,4){\circle*{3.1}}
\put(30,4){\thicklines\line(-1,-1){30}}
\put(30,4){\vector(-1,-1){15}}

\put(0,35){\thinlines\line(-1,-1){30}}
\put(0,35){\thinlines\vector(-1,-1){15}}
\put(-30,4){\circle*{3.1}}
\put(-30,4){\thinlines\line(1,-1){30}}
\put(0,-26){\thinlines\vector(-1,1){15}}

\put(0,-26){\circle*{3.1}}

\put(0,-26){\thicklines\line(0,1){30}}
\put(0,-26){\vector(0,1){15}}
\put(0,4){\circle*{3.1}}
\put(0,4){\thinlines\line(1,0){30}}
\put(30,4){\thinlines\vector(-1,0){15}}
\put(0,4){\thicklines\line(-1,0){30}}
\put(0,4){\vector(-1,0){15}}

\put(0,-26){\thinlines\line(3,1){38}}
\put(39,-13){\thinlines\vector(-3,-1){19}}
\put(0,-26){\thicklines\line(3,-1){38}}
\put(39,-39){\vector(-3,1){19}}
\put(39,-13){\circle*{3.1}}
\put(39,-39){\circle*{3.1}}
\put(39,-13){\thicklines\line(0,-1){26}}
\put(39,-13){\vector(0,-1){13}}

\put(0,-26){\thicklines\line(-3,1){38}}
\put(-39,-13){\vector(3,-1){19}}
\put(0,-26){\thinlines\line(-3,-1){38}}
\put(-39,-39){\thinlines\vector(3,1){19}}
\put(-39,-13){\circle*{3.1}}
\put(-39,-39){\circle*{3.1}}
\put(-39,-13){\thicklines\line(0,-1){26}}
\put(-39,-39){\vector(0,1){13}}
\end{picture}
\end{array}
$$
\end{Example}
\vspace{2.1cm}

\noindent {\it Notation} For each $f_i \in E({\cal{D}})\setminus
E(\tau)$ the unique cycle of the 
subgraph $\tau \cup \{f_i\}$ is denoted by $c(\tau,f_i)$. 

\begin{Proposition}\label{feb22-06} For each $f_i \in
E({\cal{D}})\setminus E(\tau)$ 
all the edges of
$c(\tau,f_i)\setminus\{f_i\}$ are oriented in the same direction and
$f_i$ is oriented in the opposite direction.
\end{Proposition}

\demo By induction on $m$, the number of subtrees $A_1,\ldots,A_m$. 
If $m=1$ the
result is easy to verify because $A_1$ is a spanning path of $G$.  
Assume $m>1$.
Consider the subgraphs
$$
\overline{G}=G\setminus\{y_1^1,\ldots,y_{i_1}^1\},\ \ 
\overline{A}_i=A_i\setminus\{y_1^1,\ldots,y_{i_1}^1\},\ i\geq 2.
$$
We set $\overline{\cal D}=(\overline{G},\overline{\cal O})$ and 
$\overline{\tau}=(\overline{A}_m,\overline{\cal O})$, where 
$\overline{\cal O}$ is the orientation induced from $\cal O$. Notice
that $\overline{G}$ is connected because $\overline{A}_m$ is a 
spanning tree of $\overline{G}$. Using the equality 
$$
V(\overline{A}_2)=\{y_{i_1+1}^2,\ldots,y_{r_1+s_1-1}^2\}=
\{z_1^1,\ldots,z_{s_1}^1,y_{i_1+2}^1,\ldots,y_{r_1}^1\}
$$
and $z_{s_1}=y_{i_1+1}^1$ it is not hard to see that 
$\overline{A}_2$ is a maximal path of $\overline{G}$ and the result 
follows by induction. Indeed  a fundamental cycle of $\tau$ is 
equal to $c(\tau,f_i)=c(\overline{\tau},f_i)$ with
$f_i\in E(\overline{\cal D})\setminus E(\overline{\tau})$ or 
$c(\tau,f_i)=c({\tau}',f_i)$ with $f_i\in E(H)\setminus
E({\tau}')$ where $H$ is the induced subgraph on 
$\{y_1^1,\ldots,y_{i_1}^1\}$ and $\tau'$ is the spanning path of 
$H$ given by $y_1^1,\ldots,y_{i_1}^1$. In the first case we apply
induction to obtain that the edges of $c(\overline{\tau},f_i)$ are 
properly oriented, in the second case it is easy to verify that 
$c({\tau}',f_i)$ has the required orientation. \QED

\begin{Theorem}\label{icd}  
$P_{\cal D}=\left(\{t_{c(\tau,f_i)}\vert f_i \in E({\cal D})\setminus
E(\tau)\}\right)$.
\end{Theorem}

\demo 
Set $E({\cal D})\setminus E(\tau)=\{f_1,\ldots,f_{q-n+1}\}$. 
Suppose without loss of generality
that $t_1,\ldots,t_{q-n+1}$ are the variables associated to 
$f_1,\ldots,f_{q-n+1}$ respectively. By Proposition~\ref{feb22-06}  
$t_{c(\tau,f_i)}=t_i-t^{\beta_i}$, where $t^{\beta_i}$ is a product 
of variables associated to edges in $\tau$. 
Let $I$ be the ideal generated by the set 
$\{t_{c(\tau,f_i)} \vert f_i \in E({\cal D})\setminus E(\tau)\}$ in
$B=k[t_1,\ldots,t_q]$.
Let $h=t^{\alpha}-t^{\beta}$ be a binomial in $P_{\cal D}$.
Thus $\overline{t_i}=\overline{t^{\beta_i}}$ in $B/I$ for $i=1,\dots,q-n+1$.
Then $\overline{h}=\overline{t^{\gamma}}-\overline{t^{\omega}}$,
where $t^{\gamma}$ and $t^{\omega}$ are products of variables
associated to 
edges of $\tau$. 
As $I\subset P_{\cal D}$, then $t^{\gamma}-t^{\omega} \in P_{\cal
D}={\rm ker}(\varphi)$. 
But $\tau$ is a tree, thus $t^{\gamma}=t^{\omega}$, and
$\overline{h}=\overline{0}$ in $B/I$. Since $P_{\cal D}$ is generated by 
binomials, $P_{\cal D}=I$. \QED

\begin{Corollary}
Assume that ${\cal D}$ is the oriented graph constructed above.
Then $P_{\cal D}$ is a homogeneous ideal 
generated by $q-n+1$ binomials corresponding to primitive cycles.
\end{Corollary}

\demo By Theorem~\ref{icd} it follows that $P_{\cal D}$ does not
contains binomials of the form $1-t^a$, i.e., $\cal D$ is acyclic. Thus
we may apply Corollary~\ref{feb23-06}. \QED

\medskip

A tournament ${\cal D}$ is a complete graph ${\cal K}_n$ with
a given orientation.

\begin{Proposition}[\rm\cite{Har}]
If ${\cal D}$ is a tournament, then ${\cal D}$ has
a spanning oriented path. 
\end{Proposition}

\begin{Proposition}[\rm\cite{ishizeki}]
If ${\cal D}$ is an acyclic tournament, then $P_{\cal D}$
is a complete intersection minimally generated by a Gr\"{o}bner
basis.
\end{Proposition}

\demo
Let $\tau$ be a spanning oriented path of ${\cal D}$, i.e., 
$\tau =\{x_1,x_2,\ldots,x_n\}$ and $(x_i,x_{i+1})$ is an edge of ${\cal D}$
for all $i<n$. Since ${\cal D}$ is acyclic, using the 
proof of Theorem~\ref{icd}, it follows
that 
$$
P_{\cal D}=(\{t_{c(\tau,f_i)} \vert f_i \in E({\cal D})\setminus E(\tau)\}),
$$
where for each $f_i \in E({\cal D})\setminus E(\tau)$, the unique
cycle of the subgraph $\tau\cup\{f_i\}$ is denoted by
$c(\tau,f_i)$. 
\QED

\medskip

Similarly we can prove the following generalization:

\begin{Proposition}
If ${\cal D}$ is an acyclic oriented graph with a spanning oriented
path, 
then $P_{\cal D}$ is a complete intersection
minimally generated by a Gr\"{o}bner basis.
\end{Proposition}

\noindent
{\bf Acknowledgments.} We thank the referees for their
careful reading of the paper and for the improvements that
they suggested.

\bibliographystyle{plain}

\end{document}